\documentclass[11pt]{article}%
\usepackage{amssymb}
\usepackage{amsmath}
\usepackage{amsfonts}
\usepackage{graphicx}
\usepackage{pgf,tikz}%
\providecommand{\U}[1]{\protect\rule{.1in}{.1in}}
\usetikzlibrary{arrows}

\newtheorem{proposition}{Proposition}[section]
\newtheorem{condition}{Condition}[section]
\newtheorem{corollary}{Corollary}[section]
\newtheorem{remark}{Remark}[section]

\begin{document}

\title{{\LARGE Logarithmic Sobolev trace inequalities}}
\author{F. Feo\thanks{Dipartimento per le Tecnologie, Universit\`a degli Studi di
Napoli ``Pathenope'', Centro Direzionale Isola C4 80143 Napoli, Italia.
E--mail: filomena.feo@uniparthenope.it} --M. R. Posteraro\thanks{Dipartimento
di Matematica e Applicazioni "R. Caccioppoli", Universit\`a degli Studi di
Napoli Federico II, Complesso Universitario Monte S. Angelo - via Cintia 80126
Napoli, Italia. E--mail: posterar@unina.it} }
\date{}
\maketitle

\begin{abstract}
{\scriptsize {\negthinspace\negthinspace\negthinspace\ }We prove a logarithmic
Sobolev trace inequality in a gaussian space and we study the trace operator
in the weighted Sobolev space }$W^{1,p}(\Omega,\gamma)$ {\scriptsize \ for
sufficiently regular domain}.{\scriptsize We exhibit examples to show the
sharpness of the results.\ Applications to PDE are also considered.}

\end{abstract}

\numberwithin{equation}{section}

\section{Introduction}

\ \ \ \ Sobolev Logarithmic inequality states that
\begin{equation}
\!\!\!\int_{%
\mathbb{R}
^{N}}\!\!\left\vert u\right\vert ^{p}\!\log\left\vert u\right\vert
d\gamma\!\leq\!\frac{p}{2}\!\int_{%
\mathbb{R}
^{N}}\!\left\vert \nabla u\right\vert ^{2}\!\left\vert u\right\vert
^{p-2}\!\text{sign }u\text{ }d\gamma\!+\!\left\Vert u\right\Vert
_{L^{p}\left(
\mathbb{R}
^{N},\gamma\right)  }^{p}\log\left\Vert u\right\Vert _{L^{p}\left(
\mathbb{R}
^{N},\gamma\right)  }\!, \label{gross con p}%
\end{equation}
where $1<p<+\infty$, $\gamma$ is the Gauss measure and $L^{p}\left(
\mathbb{R}
^{N},\gamma\right)  $ is the weighted Lebesgue space (see \S 2 for the
definitions). This inequality was first proved in \cite{[Gross]} (see also
\cite{adams articolo} for more general probability measure). It has many
applications in quantum field theory and differently from classical Sobolev
inequality it is independent of dimension and easily extends to the infinite dimension.

In terms of functional spaces inequality (\ref{gross con p}) implies the
imbedding of weighted Sobolev space $W^{1,p}(%
\mathbb{R}
^{N},\gamma)$ into the weighted Zygmund space $L^{p}\left(  \log L\right)
^{\frac{1}{2}}\left(
\mathbb{R}
^{N},\gamma\right)  .$ The imbedding holds also for $p=1$ and it is connected
with gaussian isoperimetric inequality and symmetrization (see \cite{Ledoux88}%
, \cite{eh 84} and \cite{milmann}).

\noindent For $p=+\infty$ one obtains (see \cite{[ledoux exponenzial ]} and
\cite{[aida]}) that if $u$ is a Lipschitz continuous function, then \ $u\in
L^{\infty}\left(  \log L\right)  ^{-\frac{1}{2}}\left(
\mathbb{R}
^{N},\gamma\right)  .$

\noindent This kind of imbeddings are also studied in \cite{cianchi optimal
gauss} in the more general case of rearrangement-invariant spaces.

In \cite{[mena]} a set $\Omega\subseteq\mathbb{R}^{N}$ and the space
$W_{0}^{1,p}\left(  \Omega,\gamma\right)  $ are considered using properties of
rearrangements of functions; the authors prove that if $u\in W_{0}%
^{1,p}\left(  \Omega,\gamma\right)  $ with $1\leq p<+\infty,$ then $u\in
L^{p}\left(  \log L\right)  ^{\frac{1}{2}}\left(  \Omega,\gamma\right)  $ and%
\begin{equation}
\left\Vert u\right\Vert _{L^{p}\left(  \log L\right)  ^{\frac{1}{2}}\left(
\Omega,\gamma\right)  }\leq C_{1}\left\Vert \nabla u\right\Vert _{L^{p}\left(
\Omega,\gamma\right)  }.\label{bolz Sobolevp}%
\end{equation}
Moreover if $u$ is lipschitz continuous function with$\underset
{{\scriptsize x\in\Omega,}\left\vert {\scriptsize x}\right\vert
{\scriptsize \rightarrow+\infty}}{\lim}$ $u(x)=0$ and $u_{|\partial\Omega}=0,$
then $u\!\in\!L^{\infty}\left(  \!\log L\!\right)  ^{-\frac{1}{2}}\!\left(
\Omega,\gamma\right)  $ and%
\begin{equation}
\left\Vert u\right\Vert _{L^{\infty}\left(  \log L\right)  ^{-\frac{1}{2}%
}\left(  \Omega,\gamma\right)  }\leq C_{2}\left\Vert \nabla u\right\Vert
_{L^{\infty}\left(  \Omega\right)  }.\label{bolz Sobolevinfinito}%
\end{equation}
The constants $C_{1}$,$C_{2}$ depend only on $p$ and $\gamma(\Omega).$
Analogue inequalities have been obtained in infinite dimensional case and in
the Lorentz-Zygmund spaces (see the appendix of \cite{feo hilbert+}).

A first result of our paper is to obtain (\ref{bolz Sobolevp}) when $u\in
W^{1,p}(\gamma,\Omega)$ (see \S 3); in this case, as one can expect,
smoothness assumption on $\partial\Omega$ have to be made. Besides the
continuity also the compactness of the imbedding of $W^{1,p}(\Omega,\gamma)$
in a Zygmund space is studied. As a consequence we obtain a
Poincar\'{e}-Wirtinger type inequality. We analyze also the case $p=+\infty.$
These results are sharp and counterexamples in this direction are given.
Applications of these results to PDE are also considered.

The results explained above are used to investigate Sobolev trace
ine\-qua\-lities. This kind of inequalities play a fundamental role in
problems with nonlinear boundary conditions. In the euclidean case the Sobolev
trace inequality (cf. e.g. \cite{kufner libro}) tell us that if $\Omega$ is
smooth enough and $1\leq p<N,$ then there exists a constant $C$ (depending
only on $\Omega$\ and on $p$) such that
\[
\left\Vert Tu\right\Vert _{L^{\frac{p(N-1)}{N-p}}\left(  \partial
\Omega\right)  }\leq C\left\Vert u\right\Vert _{W^{1,p}\left(  \Omega\right)
}\text{ \ \ \ \ \ for every }u\in W^{1,p}\left(  \Omega\right)  ,
\]
where $T$ is the trace operator. This kind of inequalities has been developed
via different methods and in different settings by various authors including
Besov \cite{besov intro}, Gagliardo \cite{galiardo}, Lions and Magenes
\cite{lions}. Trace inequality that involves rearrangement-invariant norms are
considered in \cite{cianchi traccia}.

To investigate about trace operator in the weighted Sobolev space
$W^{1,p}(\Omega,\gamma)$\ in \S 4 we need a Sobolev trace inequality. We prove
that if $\Omega$ is a smooth domain and $u\in C^{\infty}(\overline{\Omega})$
then
\begin{equation}%
{\displaystyle\int_{\partial\Omega}}
\left\vert u\right\vert ^{p}\log^{\frac{p}{2p^{\prime}}}(2+\left\vert
u\right\vert )\varphi dS\leq C\left\Vert u\right\Vert _{W^{1,p}(\Omega
,\gamma)}^{p}. \label{traccia lo}%
\end{equation}
This inequality is sharp and captures the spirit of the Gross inequalities:
the logarithmic function replaces the powers in this case too. We analyze also
the case $p=+\infty.$

Using (\ref{traccia lo}), we can define the trace operator and to prove
continuity and compactness of the operator $W^{1,p}(\Omega,\gamma)$ into
$L^{p}(\partial\Omega,\gamma)$ for sufficiently regular domain $\Omega
\subseteq%
\mathbb{R}
^{N}$. Moreover we prove a Poincar\'{e} trace inequality is obtained in a
suitable subspace of $W^{1,p}(\Omega,\gamma)$. We give also some applications
of these results to PDE.

An other Sobolev trace inequality is obtained in \cite{Park} as limit case of
the classical trace Sobolev inequality.

\section{\textbf{Preliminaries}}

\ \ \ \ \ In this section we recall some definitions and results which will be
useful in the following.

\subsection{Gauss measure and rearrangements}

\ \ \ \ Let $\gamma$ be the $N$-dimensional Gauss measure on $\mathbb{R}^{N}$
defined by
\[
d\gamma=\varphi\left(  x\right)  dx=\left(  2\pi\right)  ^{-\frac{N}{2}}%
\exp\left(  -\frac{\left\vert x\right\vert ^{2}}{2}\right)  \ dx,\text{
\ \ \ \ \ \ \ \ \ \ \ \ }x\in\mathbb{R}^{N}%
\]
normalized by $\gamma\left(  \mathbb{R}^{N}\right)  =1.$

\noindent We will denote by $\Phi\left(  \tau\right)  $ the Gauss measure of
the half-space $\left\{  x\in\mathbb{R}^{N}:x_{N}<\tau\right\}  :$%
\[
\!\Phi\!\left(  \!\tau\!\right)  \!=\!\gamma\!\left(  \!\left\{
\!x\!\in\!\mathbb{R}^{N}\!\!:\!x_{N}\!<\!\tau\!\right\}  \!\right)
\!=\!\frac{1}{\sqrt{2\pi}}\!\int_{-\infty}^{\tau}\!\!\!\!\exp\left(
\!-\frac{t^{2}}{2}\!\right)  dt\text{\ \ \ }\forall\tau\!\in\!\mathbb{R\cup
}\left\{  \!-\infty,\!+\infty\right\}  \!.
\]

We define\ the decreasing rearrangement with respect to Gauss measure (see
e.g. \cite{[Eh1]}) of a\ measurable function $u$ in $\Omega$ as the function
\[
u^{\circledast}\left(  s\right)  =\inf\left\{  t\geq0:\gamma_{u}\left(
t\right)  \leq s\right\}  \text{ \ \ \ \ \ \ \ \ \ \ \ \ \ \ \ \ \ \ \ \ }%
s\in\left]  0,1\right]  ,
\]
where $\gamma_{u}\left(  t\right)  =\gamma\left(  \left\{  x\in\Omega
:\left\vert u\right\vert >t\right\}  \right)  $ is the distribution function
of $u.$

\subsection{Sobolev and Zygmund space}

\ \ \ \ The weighted Lebesgue space $L^{p}(\Omega,\gamma)$ is the space of the
measurable functions $u$ on $\Omega$\ such that $\int_{\Omega}\left\vert
u\right\vert ^{p}d\gamma<+\infty.$ We recall also that the weighted Sobolev
space $W^{1,p}(\Omega,\gamma)$ for $1\leq p<+\infty$ is defined as the space
of the measurable functions $u\in L^{p}(\Omega,\gamma)$ such that there exists
$g_{1}$,...,$g_{N}\in$ $L^{p}(\Omega,\gamma)$ that verify
\[
\int_{\Omega}u\frac{\partial}{\partial x_{i}}\psi\varphi-\int_{\Omega}u\psi
x_{i}\varphi=\int_{\Omega}g_{i}\psi\varphi\text{ \ \ \ }i=1,...,N\text{
\ \ \ }\forall\psi\in D(\Omega).
\]

\noindent We stress that $u\in W^{1,p}(\Omega,\gamma)$ is a Banach space with
respect to the norm $\left\Vert u\right\Vert _{W^{1,p}(\Omega,\gamma
)}=\left\Vert u\right\Vert _{L^{p}(\Omega,\gamma)}+\left\Vert \nabla
u\right\Vert _{L^{p}(\Omega,\gamma)}.$

The Zygmund space $L^{p}(\log L)^{\alpha}(\Omega,\gamma)$ for $1\leq
p\leq+\infty$ and $\alpha\in%
\mathbb{R}
$ is the space of the measurable functions on $\Omega$ such that the quantity
\begin{equation}
||u||_{L^{p}(\log L)^{\alpha}(\Omega,\gamma)}=\!\!\left\{  \!\!%
\begin{array}
[c]{l}%
\!\left(  \!%
{\displaystyle\int_{0}^{\gamma\left(  \Omega\right)  }}
\!\left[  \!(1-\log t)^{\alpha}u^{\circledast}(t)\!\right]  ^{p}\!\!\right)
^{\frac{1}{p}}\ \ \ \ \text{if }1\leq p<+\infty\\
\!\!\underset{t\in(0,\gamma\left(  \Omega\right)  )}{\sup}\!\left[  \!(1-\log
t)^{\alpha}u^{\circledast}(t)\!\right]  \ \ \ \ \ \ \ \ \ \ \text{if\ }%
p=+\infty
\end{array}
\right.  \!\! \label{def lor-zig}%
\end{equation}
is finite. The space $L^{p}(\log L)^{\alpha}(\Omega,\gamma)$ is not trivial if
and only if $p<+\infty$ or $p=+\infty$ and $\alpha\leq0.$

\noindent The Zygmund spaces are the natural spaces in the context of Gauss
measure, because of the following property of isoperimetric function is (see
\cite{[ledoux]}):
\begin{equation}
\varphi_{1}\circ\Phi^{-1}(t)\sim t(2\log\frac{1}{t})^{\frac{1}{2}}\text{ \ for
}t\rightarrow0^{+}\text{ and }t\rightarrow1^{\text{-}}.\text{ }
\label{limite equi}%
\end{equation}

\noindent We remind same inclusion relations among Zygmund spaces. If $1\leq
r<p\leq+\infty$ and $-\infty<\alpha,\beta<+\infty$, then we get%
\[
L^{p}(\log L)^{\alpha}(\Omega,\gamma)\subseteq L^{r}(\log L)^{\beta}%
(\Omega,\gamma)\text{ .}%
\]
It is clear from definition (\ref{def lor-zig}) that the space $L^{p}(\log
L)^{\alpha}(\Omega,\gamma)$ decreases as $\alpha$ increases. For more
properties we refer to \cite{[bennet-Rudnik]}.

\subsection{Smoothness assumptions on the domain}

\ \ \ \ In this paper we deal with integrals involving the values of a
$W^{1,p}-$function on $\partial\Omega.$ To this aim we need to have a suitable
local description of the set $\Omega$ and $\partial\Omega$ is a finite union
of graphs. More precisely we will consider smooth domain $\Omega$ which
verifies the following condition (cfr. Chapter 6 of \cite{kufner libro} for
bounded domain). \newline


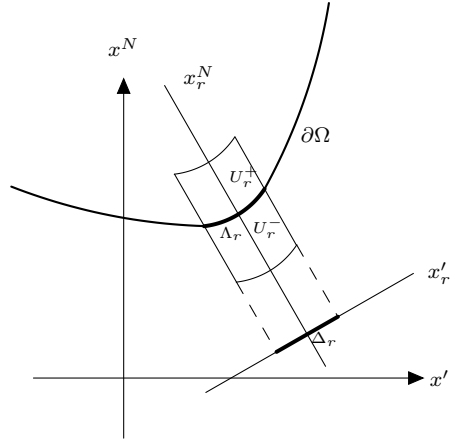
\begin{center}
\begin{tikzpicture}[line cap=round,line join=round,>=triangle 45,x=4.0mm,y=4.0mm]
\draw[->,color=black] (-3.0,0) -- (10.,0);
\draw[->,color=black] (0,-2.0) -- (0,10.0);
\draw [shift={(1.35,9.57)}] plot[domain=4.83:5.7,variable=\t]({1*2.81*cos(\t r)+0*2.81*sin(\t r)},{0*2.81*cos(\t r)+1*2.81*sin(\t r)});
\draw [shift={(3.41,5.98)}] plot[domain=4.83:5.7,variable=\t]({1*2.81*cos(\t r)+0*2.81*sin(\t r)},{0*2.81*cos(\t r)+1*2.81*sin(\t r)});
\draw [shift={(2.34,7.84)},line width=1.5 pt] plot[domain=4.83:5.7,variable=\t]({1*2.81*cos(\t r)+0*2.81*sin(\t r)},{0*2.81*cos(\t r)+1*2.81*sin(\t r)});
\draw [shift={(-12.28,15.62)},line width=0.8 pt] plot[domain=5.78:6.12,variable=\t]({1*19.36*cos(\t r)+0*19.36*sin(\t r)},{0*19.36*cos(\t r)+1*19.36*sin(\t r)});
\draw [shift={(3.30,24.40)},line width=0.8 pt] plot[domain=5.78:6.12,variable=\t]({-0.51*19.36*cos(\t r)+-0.86*19.36*sin(\t r)},{-0.86*19.36*cos(\t r)+0.51*19.36*sin(\t r)});
\draw [dash pattern=on 5pt off 5pt] (3.75,3.19)-- (5.08,0.87);
\draw [dash pattern=on 5pt off 5pt] (5.75,4.43)-- (7.12,2.04);
\draw (3.7,8.02)-- (5.74,4.42);
\draw (1.69,6.78)-- (3.75,3.19);
\draw [line width=1.5 pt] (5.08,0.87)-- (7.12,2.04);
\draw   (2.73,-0.48)-- (9.63,3.48);
\draw (1.36,9.73)-- (6.71,0.39);
\draw[color=black] (10.5,0) node {\scriptsize {$x'$}};
\draw[color=black] (10.5,3.5) node {\scriptsize {$x'_{r}$}};
\draw[color=black] (0,11) node {\scriptsize {$x^{N}$}};
\draw[color=black] (2.5,10) node {\scriptsize {$x_{r}^{N}$}};
\draw[color=black] (4.8,5) node {\tiny {$U_{r}^{-}$}};
\draw[color=black] (4,6.7) node {\tiny {$U_{r}^{+}$}};
\draw[color=black] (6.38,8.1) node {\scriptsize {$\partial\Omega$}};
\draw[color=black] (6.66,1.3) node {\tiny {$\Delta_{r}$}};
\draw[color=black] (3.58,4.8) node {\tiny {$\Lambda_{r}$}};
\end{tikzpicture}

\end{center}


\begin{condition}
\label{domain} \emph{Let\ }$\Omega$\emph{ be a domain such that there exist}

\emph{i) }$m\in%
\mathbb{N}
$\emph{ coordinate systems }$X_{r}=(x_{r}^{\prime},x_{r}^{N})$\emph{ where
}$x_{r}^{\prime}=(x_{r}^{1},...,x_{r}^{N-1})$\emph{ }

\emph{for }$r=1,2,...,m;$

\emph{ii) }$a_{i},b_{i}\in%
\mathbb{R}
\cup\left\{  \pm\infty\right\}  $\emph{ for }$i=1,...,N-1$\emph{ and }%
$m$\emph{ Lipschitz functions }$a_{r}$\emph{ in }

$\overline{\Delta_{r}}=\left\{  x_{r}^{\prime}:x_{r}^{i}\in\left(  a_{i}%
,b_{i}\right)  \text{ for }i=1,...,N-1\right\}  $\emph{ }

\emph{for }$r=1,...,N;$

\emph{iii) a number }$\beta>0$\emph{ such that the sets}%
\[
\Lambda_{r}=\left\{  (x_{r}^{\prime},x_{r}^{N})\in%
\mathbb{R}
^{N}:x_{r}^{\prime}\in\Delta_{r}\text{ and }x_{r}^{N}=a_{r}(x_{r}^{\prime
})\right\}
\]

\emph{are subsets of }$\partial\Omega$\emph{, }$\partial\Omega=\overset
{m}{\underset{r=1}{\cup}}\Lambda_{r}$\emph{ and the sets}%
\[
U_{r}^{+}=\left\{  (x_{r}^{\prime},x_{r}^{N})\in%
\mathbb{R}
^{N}:x_{r}^{\prime}\in\Delta_{r}\text{ and }a_{r}(x_{r}^{\prime})<x_{r}%
^{N}<a_{r}(x_{r}^{\prime})+\beta\right\}
\]%
\[
U_{r}^{-}=\left\{  (x_{r}^{\prime},x_{r}^{N})\in%
\mathbb{R}
^{N}:x_{r}^{\prime}\in\Delta_{r}\text{ and }a_{r}(x_{r}^{\prime})-\beta
<x_{r}^{N}<a_{r}(x_{r}^{\prime})\right\}
\]

\emph{are subset of }$\Omega$ \emph{(after a suitable orthonormal
transformation of coordinates)}$.$
\end{condition}

We observe that the set $U_{r}^{{}}=U_{r}^{+}$ $\cup$ $U_{r}^{-}$\ is an open
subset of $%
\mathbb{R}
^{N}$ and there exists an open set $U_{0}\subseteq\overline{U_{0}}%
\subset\Omega$ such that the collection $\left\{  U_{r}\right\}  _{r=0}^{m}$
is a open cover of $\Omega.$ Moreover the collection $\left\{  U_{r}\right\}
_{r=1}^{m}$ is a open cover of $\partial\Omega.$

\section{Sobolev logarithmic inequalities in $W^{1,p}(\Omega,\gamma)$}

\ \ \ \ In this section we prove continuity and compactness of imbedding of
$W^{1,p}(\Omega,\gamma)$ into $L^{p}(\log L)^{\frac{1}{2}}(\Omega,\gamma).$ We
will deal also with the case $p=+\infty.$ The first step is to obtain the
analogue of (\ref{bolz Sobolevp}) and (\ref{bolz Sobolevinfinito}) when $u\in
W^{1,p}(\Omega,\gamma)$ for $1\leq p\leq+\infty.$

\begin{proposition}
\label{prop dis sobolev w1p} (Continuity) If $u\in W^{1,p}(\Omega,\gamma)$ for
$1\leq p<+\infty$ and $\Omega$ satisfies condition \ref{domain}, then there
exists a positive constant $C$ depending only on $p$ and $\Omega$ such that
\begin{equation}
\left\Vert u\right\Vert _{L^{p}(\log L)^{\frac{1}{2}}(\Omega,\gamma)}\leq
C\left\Vert u\right\Vert _{W^{1,p}(\Omega,\gamma)}, \label{dis_w1p}%
\end{equation}
i.e. the embedding of weighted Sobolev space $W^{1,p}(\Omega,\gamma)$ into the
weighted Zygmund space $L^{p}\left(  \log L\right)  ^{\frac{1}{2}}\left(
\Omega,\gamma\right)  $ is continuous for $1\leq p<+\infty$.
\end{proposition}

To prove Proposition \ref{prop dis sobolev w1p} we need an extension operator
$P$ from $W^{1,p}(\Omega,\gamma)$ into $W^{1,p}(%
\mathbb{R}
^{N},\gamma).$ When $u\in W_{0}^{1,p}(\Omega,\gamma)$ the natural extension by
zero outside $\Omega$ is continuous without any assumptions on the regularity
of the boundary. Working with the space $W^{1,p}(\Omega,\gamma)$ the situation
is more delicate and the regularity of the boundary of $\Omega$ plays a
crucial role.

\noindent Using classical tools (see e.g. \cite{libro buttazzo} ) it is
possible to prove the existence of an extension operator $P$ from
$W^{1,p}(\Omega,\gamma)$ into $W^{1,p}(%
\mathbb{R}
^{N},\gamma)$ which is linear and continuous. The extension operator allows us
to prove the density ( for the classical case see e.g. \cite{libro buttazzo})
of $C^{\infty}(\overline{\Omega})$ in $W^{1,p}(\Omega,\gamma).$

\bigskip

\textbf{Proof of Proposition \ref{prop dis sobolev w1p}. }We consider the
extension operator $P$ from $W^{1,p}(\Omega,\gamma)$ into $W^{1,p}(%
\mathbb{R}
^{N},\gamma)$ and using (\ref{bolz Sobolevp}) we obtain for some constant $c$%
\begin{align*}
\left\Vert u\right\Vert _{L^{p}(\log L)^{\frac{1}{2}}(\Omega,\gamma)}  &  \leq
c\left\Vert Pu\right\Vert _{L^{p}(\log L)^{\frac{1}{2}}(%
\mathbb{R}
^{N},\gamma)}\leq c\left\Vert \nabla\left(  Pu\right)  \right\Vert _{L^{p}(%
\mathbb{R}
^{N},\gamma)}\\
&  \leq c\left\Vert Pu\right\Vert _{W^{1,p}(%
\mathbb{R}
^{N},\gamma)}\leq c\left\Vert u\right\Vert _{W^{1,p}(\Omega,\gamma)}%
\end{align*}
for $u\in W^{1,p}(\Omega,\gamma).$

\ \ \ \ \ \ \ \ \ \ \ \ \ \ \ \ \ \ \ \ \ \ \ \ \ \ \ \ \ \ \ \ \ \ \ \ \ \ \ \ \ \ \ \ \ \ \ \ \ \ \ \ \ \ \ \ \ \ \ \ \ \ \ \ $\ \ \ \ \ \ \ \ \ \ \ \ \ \ \ \ \ \ \ \ \ \ \ \ \ \ \ \ \ \ \ \ \ \ \ \ \ \ \ \ \ \ \ \ \ \ \ \ \ \ \ \ \ \ \ \ \ \ \ \ \ \ \ \ \ \ \ \ \ \ \ \ \ \ \ \ \ \ \ \ \ \ \square
$\bigskip

\begin{remark}
\label{remark controesempio}\emph{The space }$L^{p}(\log L)^{\frac{1}{2}%
}(\Omega,\gamma)$\emph{ obtained in the Proposition }%
\ref{prop dis sobolev w1p}\emph{ is sharp in the class of the Zygmund spaces
as the following example shows. We consider\ }$\Omega=\left\{  x\in
\mathbb{R}^{N}:x_{N}<\omega\right\}  $\emph{ with }$\omega\in%
\mathbb{R}
$ \emph{and the function }$u_{\delta}(x)=\Phi^{\delta}(x_{N})$\emph{ with}
$-\frac{1}{p}<\delta<0$\emph{. We observe that }$u_{\delta}^{\circledast
}(s)=s^{\delta}.$\emph{ By (\ref{limite equi}) we have that}%
\[
\left\Vert u_{\delta}\right\Vert _{W^{1,p}(\Omega,\gamma)}^{p}<+\infty
\Longleftrightarrow\int_{\!0}^{\!\gamma(\Omega)}\!\!\!s^{\delta p}%
(\!1\!-\!\log s\!)^{\frac{p}{2}}ds<+\infty,
\]
\emph{this means that }$u_{\delta}\in W^{1,p}(\Omega,\gamma)$ \emph{and }%
\[
\left\Vert u_{\delta}\right\Vert _{L^{p}(\log L)^{\alpha}(\Omega,\gamma
)}<+\infty\Longleftrightarrow\alpha\leq\frac{1}{2}.
\]

\end{remark}

\begin{remark}
\emph{By Proposition \ref{prop dis sobolev w1p} follows the continuity of the
embedding of Sobolev space }$W^{m,p}\!(\!\Omega,\!\!\gamma\!)$\emph{ }%
$m\!\geq1$ \emph{into the Zygmund space }$L^{\!p}\!\left(  \!\log L\!\right)
^{\!m\alpha}\left(  \!\Omega,\!\!\gamma\!\right)  $\emph{ for }$\alpha
\leq\frac{1}{2}.$ \emph{A similar result for }$\Omega=%
\mathbb{R}
^{N}$\emph{ is proved in \cite{feisner}.}
\end{remark}

Let now consider Lipschitz continuous functions.

\begin{proposition}
\label{prop infinito}If $u$ is a Lipschitz continuous function, $\Omega$
satisfies condition \ref{domain} and$\underset{{\scriptsize x\in\Omega
,}\left\vert {\scriptsize x}\right\vert {\scriptsize \rightarrow+\infty}}%
{\lim}$ $u(x)=0$, then there exists a positive constant $C$ depending only on
$\Omega$ such that
\begin{equation}
\left\Vert u\right\Vert _{L^{\infty}(\log L)^{-\frac{1}{2}}(\Omega,\gamma
)}\leq C\left(  \left\Vert \nabla u\right\Vert _{L^{\infty}(\Omega
)}+\left\Vert u\right\Vert _{L^{\infty}(\Omega)}\right)  .
\label{dis infinito}%
\end{equation}

\end{proposition}

\begin{remark}
\label{remark ottimalita infinito}\emph{The space }$L^{\infty}(\log
L)^{-\frac{1}{2}}(\Omega,\gamma)$\emph{ obtained in the Proposition
}\ref{prop infinito}\emph{ is sharp in the class of the Zygmund spaces as the
following example shows. We consider\ }$\Omega=\left\{  x\in\mathbb{R}%
^{N}:x_{N}<\omega\right\}  $\emph{ with }$\omega\in%
\mathbb{R}
$ \emph{and the function }$u_{\delta}(x)=\left(  1-\log\Phi(x_{N})\right)
^{\delta}$\emph{ with} $0<\delta\leq\frac{1}{2}$\emph{. We observe that
}$u_{\delta}^{\circledast}(s)=\left(  1-\log s\right)  ^{\delta}.$\emph{ By
(\ref{limite equi}) we have that}%
\[
\left\Vert \nabla u\right\Vert _{L^{\infty}(\Omega)}<+\infty
\Longleftrightarrow\!\!\!\underset{s\in\left(  0,\gamma(\Omega)\right)  }%
{\sup}(\!1\!-\!\log s\!)^{\delta-\frac{1}{2}}<+\infty,
\]
\emph{this means that }$u_{\delta}$ \emph{is a Lipschitz continuous function
and }%
\[
\left\Vert u_{\delta}\right\Vert _{L^{\infty}(\log L)^{\alpha}(\Omega,\gamma
)}<+\infty\Longleftrightarrow\alpha\leq-\frac{1}{2}.
\]

\end{remark}

In order to prove Proposition \ref{prop infinito} we can argue as in the proof
of Proposition \emph{\ref{prop dis sobolev w1p}:} we need the extension
operator $P$ and the inequality (\ref{bolz Sobolevinfinito}). Let us observe
that the boundary conditions $\underset{{\scriptsize x\in\Omega,}\left\vert
{\scriptsize x}\right\vert {\scriptsize \rightarrow+\infty}}{\lim}$ $u(x)=0$
and $u_{|\partial\Omega}=0$ are necessary to obtain the Polya-Sz\"{e}go
inequality for $p=+\infty,$ that is a crucial tool to prove
(\ref{bolz Sobolevinfinito}) and (\ref{dis infinito}).

\begin{proposition}
\label{proposizione compatezza}(Compactness) Let $1\leq p<+\infty$ and let
$\Omega$ satisfy condition \ref{domain}. Then the embedding of $W^{1,p}%
(\Omega,\gamma)$ into $L^{p}\log L^{\beta}(\Omega,\gamma)$ is compact if
$\beta<\frac{1}{2}.$
\end{proposition}

\textbf{Proof.} $\!$It is enough to prove the compactness of the embedding of
$W^{\!1,p}(\!\Omega,\!\!\gamma\!)$ into $L^{\substack{\!\\\!1}}\!(\!\Omega
,\!\!\gamma\!)$.$\!$ Indeed we have that any bounded set of $L^{\!p}\left(
\!\log L\!\right)  ^{\!\frac{1}{2}}\left(  \!\Omega,\!\gamma\!\right)  $ which
is precompact in $L^{1}\left(  \Omega,\gamma\right)  $ is also precompact in
$L^{p}\log L^{\beta}(\Omega,\gamma)$ with $\beta<\frac{1}{2}$\ (see e.g.
Theorem 8.23 of \cite{libro adams}).

\noindent Let be $S$ bounded set in $W^{1,p}(\Omega,\gamma),$ then $S$ is
bounded in $L^{1}(\Omega,\gamma)$ too. Using a characterization of precompact
sets of Lebesgue spaces (see e.g. Theorem 2.21 of \cite{libro adams}) we have
to prove that for any number $\varepsilon>0$ there exists a number $\delta>0$
and a subset $G\subset\subset\Omega$ such that for any $u\in S$ and every
$h\in%
\mathbb{R}
^{N}$ with $\left\vert h\right\vert <\delta$ the following conditions hold:%
\begin{equation}
a)\int_{\Omega}\left\vert \widetilde{u}(x+h)\varphi(x+h)-\widetilde
{u}(x)\varphi(x)\right\vert dx<\varepsilon\label{a}%
\end{equation}
\begin{equation}
b)\int_{\Omega-\overline{G}}\left\vert u\right\vert d\gamma<\varepsilon,
\label{b}%
\end{equation}
where $\widetilde{u}$ is the zero extension of $u$ outside $\Omega.$

\noindent Let $\varepsilon>0$ and $\Omega_{j}=\left\{  x\in\Omega
:dist(x,\partial\Omega)>\frac{1}{j}\right\}  $ for $j\in%
\mathbb{N}
.$ By (\ref{dis_w1p}) we have for some constant $c$
\begin{align*}
\int_{\Omega-\Omega_{j}}\left\vert u\right\vert d\gamma &  \leq\left(
{\displaystyle\int_{0}^{\gamma(\Omega-\Omega_{j})}}
\left[  (1-\log t)^{\frac{1}{2}}u^{\circledast}(t)\right]  ^{p}dt\right)
^{\frac{1}{p}}\left(  \int_{\Omega-\Omega_{j}}(1-\log t)^{-\frac{p^{\prime}%
}{2}}\right)  ^{\frac{1}{p^{\prime}}}\\
&  \leq c\left\Vert u\right\Vert _{W^{1,p}(\Omega,\gamma)}\left(  \int
_{\Omega-\Omega_{j}}(1-\log t)^{-\frac{p^{\prime}}{2}}\right)  ^{\frac
{1}{p^{\prime}}}%
\end{align*}
Since the Gauss measure of $\Omega$ is finite, we can choose $j$ big enough to
have
\begin{equation}
\int_{\Omega-\Omega_{j}}\left\vert u\right\vert d\gamma<\varepsilon,
\label{condizione1}%
\end{equation}
(i.e. (\ref{b}) holds) and for $h\in%
\mathbb{R}
^{N}$%
\begin{equation}
\int_{\Omega-\Omega_{j}}\left\vert \widetilde{u}(x+h)\varphi(x+h)-\widetilde
{u}(x)\varphi(x)\right\vert dx<\frac{\varepsilon}{2}.
\label{(condizione 2 parte prima)}%
\end{equation}
Let $\left\vert h\right\vert \leq\frac{1}{j},$ then $x+th\in\Omega_{2j}$ if
$x\in\Omega$ and $t\in\left[  0,1\right]  .$ Let $u\in C^{\infty}%
(\overline{\Omega}),$ we have for some constant $c$%
\begin{align}
&  \int_{\Omega_{j}}\left\vert \widetilde{u}(x+h)\varphi(x+h)-\widetilde
{u}(x)\varphi(x)\right\vert dx\label{seconda condizione parte second}\\
&  \leq\int_{\Omega_{j}}\int_{0}^{1}\left\vert \frac{d}{dt}\widetilde
{u}(x+th)\varphi(x+th)\right\vert dtdx(\text{togliere)}\nonumber\\
&  \leq\int_{\Omega_{j}}\int_{0}^{1}\left\vert \nabla\widetilde{u}%
(x+th)h\varphi(x+th)-\widetilde{u}(x+th)\varphi(x+th)\left(  x+th\right)
h\right\vert dtdx\nonumber\\
&  \leq\left\vert h\right\vert \left(  \int_{\Omega_{2j}}\left\vert
\nabla\widetilde{u}(y)\varphi(y)\right\vert dy+\int_{\Omega_{2j}}\left\vert
\widetilde{u}(y)\varphi(y)y\right\vert dy\right) \nonumber\\
&  \leq c\left\vert h\right\vert \left(  \left\Vert \nabla u\right\Vert
_{L^{p}(\Omega,\gamma)}^{p}+\left\Vert u\right\Vert _{L^{p}(\log L)^{\frac
{1}{2}}(\Omega,\gamma)}^{p}\right)  \leq c\left\vert h\right\vert \left\Vert
u\right\Vert _{W^{1,p}(\Omega,\gamma)}^{p}.\nonumber
\end{align}
In the last inequalities we have used (\ref{dis_w1p}) and the fact that
$f(x)=\left\vert x\right\vert \in L^{p^{\prime}}\left(  LogL\right)
^{-\frac{1}{2}}(\Omega,\gamma).$ Indeed since $\gamma_{f}\left(  t\right)
=1-\gamma\left(  B(0,t)\right)  ,$ one can easily check that
\[
\int_{0}^{\gamma(\Omega)}(1-\log s)^{-\frac{p^{\prime}}{2}}\left[  \left(
\left\vert x\right\vert \right)  ^{\circledast}(s)\right]  ^{p^{\prime}%
}ds=\int_{0}^{+\infty}t^{p^{\prime}}(1-\log\gamma_{f}\left(  t\right)
)^{-\frac{p^{\prime}}{2}}\gamma_{f}^{\prime}\left(  t\right)  dt<+\infty.
\]
Because of the density of $C^{\infty}(\overline{\Omega})$ in $W^{1,p}%
(\Omega,\gamma),$ (\ref{seconda condizione parte second}) holds for every $u$
in $W^{1,p}(\Omega,\gamma)$ and then for $\left\vert h\right\vert $ small
enough by (\ref{(condizione 2 parte prima)}) and
(\ref{seconda condizione parte second}) we obtain (\ref{a})

\ \ \ \ \ \ \ \ \ \ \ \ \ \ \ \ \ \ \ \ \ \ \ \ \ \ \ \ \ \ \ \ \ \ \ \ \ \ \ \ \ \ \ \ \ \ \ \ \ \ \ \ \ \ \ \ \ \ \ \ \ \ \ \ $\ \ \ \ \ \ \ \ \ \ \ \ \ \ \ \ \ \ \ \ \ \ \ \ \ \ \ \ \ \ \ \ \ \ \ \ \ \ \ \ \ \ \ \ \ \ \ \ \ \ \ \ \ \ \ \ \ \ \ \ \ \ \ \ \ \ \ \ \ \ \ \ \ \ \ \ \ \ \ \ \ \ \square
$\bigskip

\begin{remark}
\emph{Obviously the compactness results holds for }$W_{0}^{1,p}(\Omega
,\gamma)$\emph{ for any domain }$\Omega.$
\end{remark}

\begin{remark}
\emph{The compactness proved in Proposition \ref{proposizione compatezza}
implies the compact embedding of Sobolev space }$W^{m,p}\!(\!\Omega
,\!\!\gamma\!)$\emph{ }$m\!\geq1$ \emph{into the Zygmund space }%
$L^{\!p}\!\left(  \!\log L\!\right)  ^{\!m\beta}\left(  \!\Omega
,\!\!\gamma\!\right)  $\emph{ for }$\beta<\frac{1}{2}.$
\end{remark}

The compactness can be used to obtain a Poincar\'{e}-Wirtinger type inequality.

\begin{proposition}
\label{prop poincare}Let $\Omega$ be a connected domain satisfying condition
\ref{domain}. Assume $1\leq p<+\infty.$ Then there exists a positive constant
$C$, depending only on $p$ and $\Omega$, such that
\begin{equation}
\left\Vert u-u_{\Omega}\right\Vert _{L^{p}(\Omega,\gamma)}\leq C\left\Vert
\nabla u\right\Vert _{L^{p}(\Omega,\gamma)} \label{poicare}%
\end{equation}
for any $u\in W^{1,p}(\Omega,\gamma),$ where $u_{\Omega}=\frac{1}%
{\gamma(\Omega)}\int_{\Omega}ud\gamma.$
\end{proposition}

\textbf{Proof.} We precede as in the classical case. We argue by
contradiction, then there would exist for any $k\in%
\mathbb{N}
$ a function $u_{k}\in W^{1,p}(\Omega,\gamma)$ such that%
\[
\left\Vert u_{k}-\left(  u_{k}\right)  _{\Omega}\right\Vert _{L^{p}%
(\Omega,\gamma)}>k\left\Vert \nabla u_{k}\right\Vert _{L^{p}(\Omega,\gamma)}.
\]
We renormalize by defining
\begin{equation}
v_{k}=\frac{u_{k}-\left(  u_{k}\right)  _{\Omega}}{\left\Vert u_{k}-\left(
u_{k}\right)  _{\Omega}\right\Vert _{L^{p}(\Omega,\gamma)}}. \label{vk}%
\end{equation}
Then%
\[
\left(  v_{k}\right)  _{\Omega}=0,\left\Vert v_{k}\right\Vert _{L^{p}%
(\Omega,\gamma)}=1
\]
and
\begin{equation}
\text{and }\left\Vert \nabla v_{k}\right\Vert _{L^{p}(\Omega,\gamma)}<\frac
{1}{k}. \label{limitatezza}%
\end{equation}
In particular the functions $\left\{  v_{k}\right\}  _{k\in%
\mathbb{N}
}$ are bounded in $W^{1,p}(\Omega,\gamma).$ Then by the previous theorem there
exists a subsequence still denoted by $\left\{  v_{k}\right\}  _{k\in%
\mathbb{N}
}$ and a function $v$ such that
\[
v_{k}\rightarrow v\text{ \ \ \ in }L^{p}(\Omega,\gamma)\text{.}%
\]
Moreover by (\ref{vk}) it follows that
\begin{equation}
v_{\Omega}=0\text{ and }\left\Vert v\right\Vert _{L^{p}(\Omega,\gamma)}=1.
\label{valore medio}%
\end{equation}
On the other hand, (\ref{limitatezza}) implies for any $\psi\in C_{0}^{\infty
}(\Omega)$ and $i=1,...,N$%
\begin{align*}
\int_{\Omega}v\frac{\partial\psi}{\partial x_{i}}\varphi dx-\int_{\Omega}v\psi
x_{i}\varphi dx  &  =\underset{k\rightarrow+\infty}{\lim}\left(  \int_{\Omega
}v_{k}\frac{\partial\psi}{\partial x_{i}}\varphi dx-\int_{\Omega}v_{k}\psi
x_{i}\varphi dx\right) \\
&  =\underset{k\rightarrow+\infty}{\lim}-\int_{\Omega}\frac{\partial v_{k}%
}{\partial x_{i}}\psi\varphi dx=0.
\end{align*}
Consequently $v\in W^{1,p}(\Omega,\gamma)$ and $\nabla v=0$ a.e. Then $v$ is
constant since $\Omega$ is connected. In particular by the first estimate in
(\ref{valore medio}) we must have $v\equiv0;$ in which case $\left\Vert
v\right\Vert _{L^{p}(\Omega,\gamma)}=0.$ This contradiction establishes the
estimate (\ref{poicare}).

\ \ \ \ \ \ \ \ \ \ \ \ \ \ \ \ \ \ \ \ \ \ \ \ \ \ \ \ \ \ \ \ \ \ \ \ \ \ \ \ \ \ \ \ \ \ \ \ \ \ \ \ \ \ \ \ \ \ \ \ \ \ \ \ $\ \ \ \ \ \ \ \ \ \ \ \ \ \ \ \ \ \ \ \ \ \ \ \ \ \ \ \ \ \ \ \ \ \ \ \ \ \ \ \ \ \ \ \ \ \ \ \ \ \ \ \ \ \ \ \ \ \ \ \ \ \ \ \ \ \ \ \ \ \ \ \ \ \ \ \ \ \ \ \ \ \ \square
$\bigskip

\begin{remark}
\label{remark sottospazi}\emph{The previous proof works in a more general
case. Let }$\Omega$\emph{ be a connected domain satisfying condition
\ref{domain} and let }$V\subset W^{1,p}(\Omega,\gamma)$\emph{ be a linear
subspace of }$W^{1,p}(\Omega,\gamma)$\emph{ with }$1\leq p<+\infty$\emph{
which is closed and such that the only constant function belonging to }%
$V$\emph{ is the function which is identically zero. Then there exists a
positive constant }$C$\emph{, depending only on }$p$\emph{ and }$\Omega
$\emph{, such that }%
\[
\left\Vert v\right\Vert _{L^{p}(\Omega,\gamma)}\leq C\left(  \int_{\Omega
}\overset{N}{\underset{i=1}{%
{\displaystyle\sum}
}}\left\vert \frac{\partial v}{\partial x_{i}}\right\vert ^{p}d\gamma\right)
^{\frac{1}{p}}\text{ \ \ \ \ }\forall v\in V.
\]

\end{remark}

\begin{remark}
\emph{(}Application to PDE\emph{) Let }$\Omega$\emph{ be a connected domain
satisfying condition \ref{domain}. Let us consider the semicoercive
homogeneous Neumann problem }%
\begin{equation}
\left\{
\begin{array}
[c]{ll}%
-\left(  u_{x_{i}}\varphi\right)  _{x_{i}}=f\varphi & \text{\emph{\mbox{in}}
}\Omega\\
& \\
\frac{\partial u}{\partial\nu}=0\quad & \text{\emph{\mbox{on}} }\partial
\Omega,
\end{array}
\right.  \label{p1}%
\end{equation}
\emph{where }$f\in L^{2}(\log L)^{-\frac{1}{2}}(\Omega,\gamma)$ \emph{and
}$\nu$\emph{ is the external normal. Using classical tools (see e.g.
\cite{libro buttazzo} Theorem 6.2.3) and inequalities (\ref{dis_w1p}) and
(\ref{poicare}) it follows that problem (\ref{p1}) has a weak solution in
}$W^{1,2}(\Omega,\gamma)$\emph{ if and only if }$\int_{\Omega}fd\gamma
=0.$\emph{ In particular there exists a unique weak solution in }$X=\left\{
u\in W^{1,2}(\Omega,\gamma):\int_{\Omega}ud\gamma=0\right\}  $\emph{ by
Lax-Milgram theorem}$.$

\emph{We consider also the following eigenvalue problem related to the
equation of quantum harmonic oscillator}%
\begin{equation}
\left\{
\begin{array}
[c]{ll}%
-\left(  u_{x_{i}}\varphi\right)  _{x_{i}}=\lambda u & \text{\emph{\mbox{in}}
}\Omega\\
& \\
\frac{\partial u}{\partial\nu}=0\quad & \text{\emph{\mbox{on}} }\partial
\Omega.
\end{array}
\right.  \label{autovalori newmann}%
\end{equation}
\emph{Arguing in a classical way (see e.g. \cite{libro buttazzo} Theorem
8.6.1), using inequality (\ref{poicare}) and the compactness of the embedding
from }$W^{1,2}(\Omega,\gamma)$\emph{ into }$L^{2}(\Omega,\gamma),$\emph{ it
follows that} \emph{there exists an increasing sequence of eigenvalues of the
problem (\ref{autovalori newmann}) which tends to infinity and a Hilbertian
basis of eigenfunctions in }$L^{2}(\Omega,\gamma)$\emph{. Moreover for
}$\lambda_{1}=0,$\emph{ the corresponding eigenfunction }$u_{1}=const\neq
0$\emph{ and the first nontrivial eigenvalue }$\lambda_{2}$\emph{\ has the
following characterization}%
\[
\lambda_{2}=\min\left\{  \frac{\left\Vert \nabla u\right\Vert _{L^{2}%
(\Omega,\gamma)}}{\left\Vert u\right\Vert _{L^{2}(\Omega,\gamma)}},u\in
W^{1,2}(\Omega,\gamma):\int_{\Omega}ud\gamma=0\right\}  .
\]

\end{remark}

\section{Sobolev logarithmic trace inequalities}

\ \ \ \ In this section we deal with integrals involving the values of a
$C^{\infty}-$function on $\partial\Omega.$ We prove that a certain integral of
the function on $\partial\Omega$ is bounded by the $W^{1,p}-$norm on $\Omega.$
This inequality will be crucial to define trace operator (see \S 5).

\begin{proposition}
\label{prop dis lplogl} Let $\Omega$ be a domain satisfying condition\emph{
\ref{domain}} and $1\leq p<+\infty.$ For every $u\in C^{\infty}(\overline
{\Omega})$ there exists a positive constant $C$ depending only on $p$ and
$\Omega$ such that%
\begin{equation}%
{\displaystyle\int_{\partial\Omega}}
\left\vert u\right\vert ^{p}\log^{\frac{p-1}{2}}(2+\left\vert u\right\vert
)\varphi dS\leq C\left\Vert u\right\Vert _{W^{1,p}(\Omega,\gamma)}^{p}
\label{dis lplogl}%
\end{equation}

\end{proposition}

\begin{remark}
\emph{We obtain the same result if we replace the first member of
(\ref{dis lplogl}) with the quantity }$%
{\displaystyle\int_{\partial\Omega}}
u^{p}\left(  \log^{+}(\left\vert u\right\vert )\right)  ^{\frac{p-1}{2}%
}\varphi dS.$
\end{remark}

\textbf{Proof. }Following classical tools (see Chapter 6 of \cite{kufner
libro} ) it is enough to prove the existence of a constant $C_{T}>0$ such that
for any function $u\in C^{\infty}(\overline{\Omega})$ whose supports is in
$\Lambda_{r}\cup U_{r}^{+}$ we have (\ref{dis lplogl}). After suitable
transformation that maps $\Delta_{r}\times\left]  0,\beta\right[  \ $onto
$U_{r}^{+}$ and $\Delta_{r}\times\left\{  0\right\}  \ $onto $\Lambda_{r},$ we
can reduce to consider $u$ such that the support is in $\Delta_{r}%
\times\left[  0,\beta\right[  $. Then it is sufficient to prove the existence
of a constant $C>0$ such that for any function $u\in C^{\infty}(\overline
{\Delta_{r}}\times\left[  0,\beta\right[  )$ whose supports is in $\Delta
_{r}\times\left[  0,\beta\right[  $%

\begin{equation}%
{\displaystyle\int_{\Delta_{r}}}
\left\vert u(x_{r}^{\prime},0)\right\vert ^{p}\log^{\frac{p-1}{2}%
}(2+\left\vert u(x_{r}^{\prime},0)\right\vert )\varphi(x_{r}^{\prime},0)\text{
}dx_{r}^{\prime}\leq C\left\Vert u\right\Vert _{W^{1,p}(\Delta_{r}%
\times\left]  0,\beta\right[  ,\gamma)}^{p}. \label{dis traccia locale}%
\end{equation}
holds. In (\ref{dis traccia locale}) we have denoted by $u$ the composition of
$u$ with the change of coordinates.

Now we prove (\ref{dis traccia locale})\textbf{.} For some constant $c$ that
can varies from line to line we have
\begin{align}
&
{\displaystyle\int_{\Delta_{r}}}
\left\vert u(x_{r}^{\prime},0)\right\vert ^{p}\log^{\frac{p-1}{2}%
}(2+\left\vert u(x_{r}^{\prime},0)\right\vert )\text{ }\varphi(x_{r}^{\prime
},0)\text{ }dx_{r}^{\prime}\label{prima}\\
&  \leq c(A_{1}+A_{2}+A_{3})\nonumber
\end{align}
where
\begin{align*}
A_{1}  &  =%
{\displaystyle\int_{\Delta_{r}}}
{\displaystyle\int_{\beta}^{0}}
p\left\vert u(x_{r}^{\prime},x_{r}^{N})\right\vert ^{p-1}\log^{\frac{p-1}{2}%
}(2+\left\vert u(x_{r}^{\prime},x_{r}^{N})\right\vert )\left\vert
\frac{\partial u}{\partial x_{r}^{N}}(x_{r}^{\prime},x_{r}^{N})\right\vert
\varphi(x_{r}^{\prime},x_{r}^{N})\text{ }dx_{r}^{N}\text{ }dx_{r}^{\prime}\\
A_{2}  &  =%
{\displaystyle\int_{\Delta_{r}}}
{\displaystyle\int_{\beta}^{0}}
\frac{p-1}{2}\left\vert u(x_{r}^{\prime},x_{r}^{N})\right\vert ^{p}\frac
{\log^{\frac{p-1}{2}-1}(2+\left\vert u(x_{r}^{\prime},x_{r}^{N})\right\vert
)}{2+\left\vert u(x_{r}^{\prime},x_{r}^{N})\right\vert }\left\vert
\frac{\partial u}{\partial x_{r}^{N}}(x_{r}^{\prime},x_{r}^{N})\right\vert
\varphi(x_{r}^{\prime},x_{r}^{N})\text{ }dx_{r}^{N}\text{ }dx_{r}^{\prime}\\
A_{3}  &  =%
{\displaystyle\int_{\Delta_{r}}}
{\displaystyle\int_{\beta}^{0}}
\left\vert u(x_{r}^{\prime},x_{r}^{N})\right\vert ^{p}\log^{\frac{p-1}{2}%
}(2+\left\vert u(x_{r}^{\prime},x_{r}^{N})\right\vert )\varphi(x_{r}^{\prime
},x_{r}^{N})\text{ }\left\vert x_{r}^{N}\right\vert \text{ }dx_{r}^{N}\text{
}dx_{r}^{\prime}.
\end{align*}
We observe that the function $f(x)=x_{r}^{N}\in L^{\infty}\left(  \log
L\right)  ^{-\frac{1}{2}}(\Delta_{r}\times\left]  0,\beta\right[  ,\gamma).$
Indeed $\gamma_{f}\left(  t\right)  =2\Phi\left(  -t\right)  $ and using
(\ref{limite equi}) we have
\begin{align*}
\underset{t\in(0,\gamma(\Delta_{r}\times\left]  0,\beta\right[  ))}{\sup
}(1-\log t)^{-\frac{1}{2}}f^{\circledast}(t)  &  =\underset{t\in
(0,\gamma(\Delta_{r}\times\left]  0,\beta\right[  ))}{\sup}(1-\log
t)^{-\frac{1}{2}}\left(  -\Phi^{-1}\left(  \frac{t}{2}\right)  \right) \\
&  \leq c\underset{t\in(0,\gamma(\Delta_{r}\times\left]  0,\beta\right[
))}{\sup}(1-\log t)^{-\frac{1}{2}}(2\log\frac{2}{t})^{\frac{1}{2}}<+\infty.
\end{align*}
Then we obtain
\begin{equation}
A_{3}\leq c\left\Vert u\right\Vert _{L^{p}\left(  \log L\right)  ^{\frac
{1}{2p^{\prime}}}(\Delta_{r}\times\left]  0,\beta\right[  ,\gamma)}%
^{p}\left\Vert x_{r}^{N}\right\Vert _{L^{\infty}\left(  \log L\right)
^{-\frac{1}{2}}(\Delta_{r}\times\left]  0,\beta\right[  ,\gamma)}. \label{A3}%
\end{equation}
Moreover using H\"{o}lder inequality, we obtain
\begin{align}
A_{1}  &  \leq c\left(
{\displaystyle\int_{\Delta_{r}\times\left]  0,\beta\right[  }}
\left\vert u(x_{r}^{\prime},x_{r}^{N})\right\vert ^{p}\log^{\frac{p}{2}%
}(2+\left\vert u(x_{r}^{\prime},x_{r}^{N})\right\vert )\varphi(x_{r}^{\prime
},x_{r}^{N})\text{ }dx_{r}^{N}\text{ }dx_{r}^{\prime}\right)  ^{\frac
{1}{p^{\prime}}}\label{A1}\\
&  \times\left(
{\displaystyle\int_{\Omega}}
\left\vert \frac{\partial u}{\partial x_{r}^{N}}(x_{r}^{\prime},x_{r}%
^{N})\right\vert ^{p}\varphi(x_{r}^{\prime},x_{r}^{N})\text{ }dx_{r}^{N}\text{
}dx_{r}^{\prime}\right)  ^{\frac{1}{p}}\nonumber\\
A_{2}  &  \leq c\left(
{\displaystyle\int_{\Delta_{r}\times\left]  0,\beta\right[  }}
\left\vert u(x_{r}^{\prime},x_{r}^{N})\right\vert ^{p}\log^{\left(  \frac
{p-1}{2}-1\right)  p^{\prime}}(2+\left\vert u(x_{r}^{\prime},x_{r}%
^{N})\right\vert )\varphi(x_{r}^{\prime},x_{r}^{N})\text{ }dx_{r}^{N}\text{
}dx_{r}^{\prime}\right)  ^{\frac{1}{p^{\prime}}}\times\label{A2}\\
&  \times\left(
{\displaystyle\int_{\Omega}}
\left\vert \frac{\partial u}{\partial x_{r}^{N}}(x_{r}^{\prime},x_{r}%
^{N})\right\vert ^{p}\varphi(x_{r}^{\prime},x_{r}^{N})\text{ }dx_{r}^{N}\text{
}dx_{r}^{\prime}\right)  ^{\frac{1}{p}}\nonumber
\end{align}
We observe that%
\begin{align*}
&
{\displaystyle\int_{\Delta_{r}\times\left]  0,\beta\right[  }}
\left\vert u(x_{r}^{\prime},x_{r}^{N})\right\vert ^{p}\log^{\left(  \frac
{p-1}{2}-1\right)  p^{\prime}}(2+\left\vert u(x_{r}^{\prime},x_{r}%
^{N})\right\vert )\varphi(x_{r}^{\prime},x_{r}^{N})\text{ }dx_{r}^{N}\text{
}dx_{r}^{\prime}\\
&  \leq c%
{\displaystyle\int_{\Delta_{r}\times\left]  0,\beta\right[  }}
\left\vert u(x_{r}^{\prime},x_{r}^{N})\right\vert ^{p}\log^{\frac{p}{2}%
}(2+\left\vert u(x_{r}^{\prime},x_{r}^{N})\right\vert )\varphi(x_{r}^{\prime
},x_{r}^{N})\text{ }dx_{r}^{N}\text{ }dx_{r}^{\prime}%
\end{align*}
and
\begin{equation}
\left(
{\displaystyle\int_{\Delta_{r}\times\left]  0,\beta\right[  }}
\left\vert u(x_{r}^{\prime},x_{r}^{N})\right\vert ^{p}\log^{\frac{p}{2}%
}(2+\left\vert u(x_{r}^{\prime},x_{r}^{N})\right\vert )\varphi(x_{r}^{\prime
},x_{r}^{N})\text{ }dx_{r}^{N}\text{ }dx_{r}^{\prime}\right)  \leq\label{aa}%
\end{equation}%
\begin{align*}
&  =%
{\displaystyle\int_{0}^{\gamma(\Delta_{r}\times\left]  0,\beta\right[  )}}
\left[  u^{\circledast}(t)\log^{\frac{1}{2}}(2+u^{\circledast}(t))\right]
^{p}dt\\
&  \leq c\left(
{\displaystyle\int_{0}^{\gamma(\Delta_{r}\times\left]  0,\beta\right[  )}}
\left[  (1-\log t)^{\frac{1}{2}}u^{\circledast}(t)\right]  ^{p}dt\right)  ,
\end{align*}
because $\log(2+u^{\circledast}(t))$ is dominated by a multiple of $(1-\log
t)$. Indeed $L^{p}\left(  \log L\right)  ^{\frac{1}{2}}\subset L^{p}\subset
L^{p,\infty}$, then $u^{\circledast}(t)\leq ct^{-\frac{1}{p}}$ for some
positive constant$.$

Putting (\ref{A3})-(\ref{aa}) in (\ref{prima}) and using Proposition
\ref{prop dis sobolev w1p} we have
\begin{align*}
&
{\displaystyle\int_{\Delta_{r}}}
\left\vert u(x_{r}^{\prime},0)\right\vert ^{p}\log^{\frac{p-1}{2}%
}(2+\left\vert u(x_{r}^{\prime},0)\right\vert )\varphi(x_{r}^{\prime}%
,0)dx_{r}^{\prime}\\
\!\!  &  \leq\!c\!\left\Vert u\right\Vert \!_{\!L^{p}\!\left(  \log L\right)
^{\frac{1}{2}}\!(\Delta_{r}\times\left]  0,\beta\right[  ,\gamma)}%
^{p-1}\!\!\!\!\left\Vert \nabla u\right\Vert \!\!_{L^{p}(\Delta_{r}%
\times\left]  0,\beta\right[  ,\gamma)}+\!c\left\Vert u\right\Vert
\!_{\!L^{p}\left(  \log L\right)  ^{\frac{1}{2p^{\prime}}\!}(\Delta_{r}%
\times\left]  0,\beta\right[  ,\gamma)\!}^{p}\!\!\!\!\!\left\Vert x_{r}%
^{N}\right\Vert \!_{\!L^{\infty}\!\left(  \log L\right)  ^{-\frac{1}{2}%
}\!(\Delta_{r}\times\left]  0,\beta\right[  ,\gamma)}\\
\!\!  &  \leq\!c\left\Vert u\right\Vert _{W^{1,p}(\Delta_{r}\times\left]
0,\beta\right[  ,\gamma)}^{p-1}\left\Vert \nabla u\right\Vert _{L^{p}%
(\Delta_{r}\times\left]  0,\beta\right[  ,\gamma)}+c\left\Vert u\right\Vert
_{W^{1,p}(\Delta_{r}\times\left]  0,\beta\right[  ,\gamma)}^{p}\\
&  \leq c\left\Vert u\right\Vert _{W^{1,p}(\Delta_{r}\times\left]
0,\beta\right[  ,\gamma)}^{p}.
\end{align*}

\ \ \ \ \ \ \ \ \ \ \ \ \ \ \ \ \ \ \ \ \ \ \ \ \ \ \ \ \ \ \ \ \ \ \ \ \ \ \ \ \ \ \ \ \ \ \ \ \ \ \ \ \ \ \ \ \ \ \ \ \ \ \ \ $\ \ \ \ \ \ \ \ \ \ \ \ \ \ \ \ \ \ \ \ \ \ \ \ \ \ \ \ \ \ \ \ \ \ \ \ \ \ \ \ \ \ \ \ \ \ \ \ \ \ \ \ \ \ \ \ \ \ \ \ \ \ \ \ \ \ \ \ \ \ \ \ \ \ \ \ \ \ \ \ \ \ \square
$\bigskip

\begin{remark}
\label{remark ottomalita traccia}\emph{In (\ref{dis lplogl})the exponent
}$\frac{p-1}{2}$ \emph{of the logarithmic is sharp as the following example
shows. We consider\ }$\Omega=\left\{  x\in\mathbb{R}^{N}:x_{N}<\omega\right\}
$\emph{ with }$\omega\in%
\mathbb{R}
$ \emph{and the function }$u_{\delta}(x)=\Phi^{\delta}(x_{N})$\emph{ with}
$-\frac{1}{p}<\delta<0$ \emph{as in Remark \ref{remark controesempio}. We have
that }%
\[%
{\displaystyle\int_{\partial\Omega}}
\left\vert u\right\vert ^{p}\log^{\beta}(2+\left\vert u\right\vert )\varphi
dS=A_{1}+A_{2}+A_{3}%
\]
\emph{where}
\begin{align*}
A_{1}  &  =%
{\displaystyle\int}
{\displaystyle\int_{\Omega}}
p\left\vert u(x^{^{\prime}},x_{N})\right\vert ^{p-1}\log^{\beta}(2+\left\vert
u(x^{\prime},x_{N})\right\vert )\left\vert \frac{\partial u}{\partial x_{N}%
}(x_{r}^{\prime},x_{N})\right\vert \varphi(x^{\prime},x_{N})dx^{\prime}%
dx_{N}\\
A_{2}  &  =%
{\displaystyle\int}
{\displaystyle\int_{\Omega}}
\frac{p-1}{2}\left\vert u(x^{\prime},x_{N})\right\vert ^{p}\frac{\log
^{\beta-1}(2+\left\vert u(x^{\prime},x_{N})\right\vert )}{2+\left\vert
u(x^{\prime},x_{N})\right\vert }\left\vert \frac{\partial u}{\partial x_{N}%
}(x^{\prime},x_{N})\right\vert \varphi(x^{\prime},x_{N})dx^{\prime}dx_{N}\\
A_{3}  &  =-%
{\displaystyle\int}
{\displaystyle\int_{\Omega}}
\left\vert u(x_{r}^{\prime},x_{N})\right\vert ^{p}\log^{\frac{p-1}{2}%
}(2+\left\vert u(x^{\prime},x_{N})\right\vert )\varphi(x^{\prime},x_{N})\text{
}x_{N}\text{ }dx^{\prime}dx_{N}.
\end{align*}
\emph{By (\ref{limite equi}) we have that}%
\[
A_{1}<+\infty\Longleftrightarrow\beta\leq\frac{p-1}{2}<\frac{p}{2}.
\]
\emph{Moreover }$A_{2}\leq cA_{1}$\emph{ and }%
\[
A_{3}<+\infty\Longleftrightarrow\beta\leq\frac{p}{2},
\]
\emph{because }$x_{N}\in L^{\infty}\left(  \log L\right)  ^{-\frac{1}{2}%
}(\Omega,\gamma).$
\end{remark}

If $p=+\infty$ we can prove the following result.

\begin{proposition}
\label{prop traccia infinito}Let $\Omega$ be a domain satisfying
condition\emph{ \ref{domain}}. For every $u\in C^{\infty}(\overline{\Omega})$
such that $\underset{x\in\Omega,\left\vert x\right\vert \rightarrow+\infty
}{\lim}u=0$ and every $\lambda\in\left(  0,1\right)  $, there exists a
positive constant depending on $\Omega$ and $\lambda$ such that
\begin{align}%
{\displaystyle\int_{\partial\Omega}}
\exp\left(  \lambda\left\vert u\right\vert ^{2}\right)  \varphi\text{ }dS  &
\leq C\exp\left[  \left(  \left\Vert \nabla u\right\Vert _{L^{\infty}(\Omega
)}+\left\Vert u\right\Vert _{L^{\infty}(\Omega)}\right)  ^{2}\right]
\label{traccia infinito}\\
&  \left(  \left\Vert \nabla u\right\Vert _{L^{\infty}(\Omega)}\left(
\left\Vert \nabla u\right\Vert _{L^{\infty}(\Omega)}+\left\Vert u\right\Vert
_{L^{\infty}(\Omega)}\right)  +1\right)  .\nonumber
\end{align}

\end{proposition}

\textbf{Proof of Proposition \ref{prop traccia infinito}. }As in the proof of
Proposition \ref{prop dis lplogl} it is sufficient to prove for any functions
$u\in C^{\infty}(\overline{\Delta_{r}}\times\left[  0,\beta\right[  )$ whose
supports is in $\Delta_{r}\times\left[  0,\beta\right[  $ and any $\lambda
\in\left(  0,1\right)  $ the following inequality%

\begin{equation}%
{\displaystyle\int_{\Delta_{r}}}
\exp\left(  \lambda\left\vert u(x_{r}^{\prime},0)\right\vert ^{2}\right)
\varphi(x_{r}^{\prime},0)dx_{r}^{\prime}\leq C\exp\left[  \left(  \left\Vert
\nabla u\right\Vert _{L^{\infty}(\overline{\Delta_{r}}\times\left[
0,\beta\right[  )}+\left\Vert u\right\Vert _{L^{\infty}(\overline{\Delta_{r}%
}\times\left[  0,\beta\right[  )}\right)  ^{2}\right]  \times
\label{infinito locale}%
\end{equation}%
\[
\times\left(  \left\Vert \nabla u\right\Vert _{L^{\infty}(\overline{\Delta
_{r}}\times\left[  0,\beta\right[  )}\left(  \left\Vert \nabla u\right\Vert
_{L^{\infty}(\overline{\Delta_{r}}\times\left[  0,\beta\right[  )}+\left\Vert
u\right\Vert _{L^{\infty}(\overline{\Delta_{r}}\times\left[  0,\beta\right[
)}\right)  +1\right)
\]
holds for some positive constant $C$ not depending on $u$ .

Now we prove (\ref{infinito locale}). For some constant $c$ that can varies
from line to line we have
\begin{equation}%
{\displaystyle\int_{\Delta_{r}}}
\exp\left(  \lambda\left\vert u(x_{r}^{\prime},0)\right\vert ^{2}\right)
\varphi(x_{r}^{\prime},0)dx_{r}^{\prime}\leq c(B_{1}+B_{2}) \label{B1_B2}%
\end{equation}
where
\begin{align*}
B_{1}  &  =%
{\displaystyle\int_{\Delta_{r}}}
{\displaystyle\int_{\beta}^{0}}
\lambda\left\vert u(x_{r}^{\prime},0)\right\vert \exp\left(  \lambda\left\vert
u(x_{r}^{\prime},0)\right\vert ^{2}\right)  \left\vert \frac{\partial
u}{\partial x_{r}^{N}}(x_{r}^{\prime},x_{r}^{N})\right\vert \varphi
(x_{r}^{\prime},x_{r}^{N})\text{ }dx_{r}^{N}\text{ }dx_{r}^{\prime}\\
B_{2}  &  =%
{\displaystyle\int_{\Delta_{r}}}
{\displaystyle\int_{\beta}^{0}}
\exp\left(  \lambda\left\vert u(x_{r}^{\prime},0)\right\vert ^{2}\right)
\varphi(x_{r}^{\prime},x_{r}^{N})\text{ }\left\vert x_{r}^{N}\right\vert
\text{ }dx_{r}^{N}\text{ }dx_{r}^{\prime}.
\end{align*}
Moreover, since $u^{\circledast}(t)\leq\left\Vert u\right\Vert _{L^{\infty
}(\log L)^{-\frac{1}{2}}(\overline{\Delta_{r}}\times\left[  0,\beta\right[
)}(1-\log t)^{\frac{1}{2}}$ in $\overline{\Delta_{r}}\times\left[
0,\beta\right[  $, we have%
\begin{align}
B_{1}  &  \leq\lambda\left\Vert \nabla u\right\Vert _{L^{\infty}%
(\overline{\Delta_{r}}\times\left[  0,\beta\right[  )}%
{\displaystyle\int_{\Delta_{r}}}
{\displaystyle\int_{\beta}^{0}}
\left\vert u(x_{r}^{\prime},0)\right\vert \exp\left(  \lambda\left\vert
u(x_{r}^{\prime},0)\right\vert ^{2}\right)  \varphi(x_{r}^{\prime},x_{r}%
^{N})\text{ }dx_{r}^{N}\text{ }dx_{r}^{\prime}\label{b1}\\
&  =\lambda\left\Vert \nabla u\right\Vert _{L^{\infty}(\overline{\Delta_{r}%
}\times\left[  0,\beta\right[  )}%
{\displaystyle\int_{0}^{\gamma(\overline{\Delta_{r}}\times\left[
0,\beta\right[  )}}
\left\vert u^{\circledast}(t)\right\vert \exp\left(  \lambda\left\vert
u^{\circledast}(t)\right\vert ^{2}\right)  dt\nonumber\\
&  \leq\lambda\exp\left(  \left\Vert u\right\Vert _{L^{\infty}(\log
L)^{-\frac{1}{2}}(\overline{\Delta_{r}}\times\left[  0,\beta\right[  )}%
^{2}\right)  \left\Vert \nabla u\right\Vert _{L^{\infty}(\overline{\Delta_{r}%
}\times\left[  0,\beta\right[  )}\left\Vert u\right\Vert _{L^{\infty}(\log
L)^{-\frac{1}{2}}(\overline{\Delta_{r}}\times\left[  0,\beta\right[  )}%
\times\nonumber\\
&  \times%
{\displaystyle\int_{0}^{\gamma(\overline{\Delta_{r}}\times\left[
0,\beta\right[  )}}
\frac{(1-\log t)^{\frac{1}{2}}}{t^{\lambda}}dt.\nonumber
\end{align}
and
\begin{equation}
\!\!B_{2}\!\leq\!\left\Vert x_{r}^{N}\right\Vert \!_{\!L^{\infty}\!\left(
\log L\right)  ^{-\frac{1}{2}}\!(\Delta_{r}\times\left]  0,\beta\right[
,\gamma)}\!\exp\left(  \!\left\Vert u\right\Vert \!_{\!L^{\infty}(\log
L)^{-\frac{1}{2}}(\overline{\Delta_{r}}\times\left[  0,\beta\right[  )}%
^{2}\!\right)  \!%
{\displaystyle\int_{0}^{\gamma(\overline{\Delta_{r}}\times\left[
0,\beta\right[  )}}
\!\!\frac{(1-\log t)^{\frac{1}{2}}}{t^{\lambda}}dt. \label{b2}%
\end{equation}
For any $\lambda\in\left(  0,1\right)  $ the integrals in (\ref{b1}) and
(\ref{b2}) are finite and%
\begin{equation}
B_{1}\leq c\left\Vert \nabla u\right\Vert _{L^{\infty}(\overline{\Delta_{r}%
}\times\left[  0,\beta\right[  )}\left\Vert u\right\Vert _{L^{\infty}(\log
L)^{-\frac{1}{2}}(\overline{\Delta_{r}}\times\left[  0,\beta\right[  )}%
\exp\left(  \left\Vert u\right\Vert _{L^{\infty}(\log L)^{-\frac{1}{2}%
}(\overline{\Delta_{r}}\times\left[  0,\beta\right[  )}^{2}\right)
\label{stima b1}%
\end{equation}
and
\begin{equation}
B_{2}\leq c\exp\left(  \left\Vert u\right\Vert _{L^{\infty}(\log L)^{-\frac
{1}{2}}(\overline{\Delta_{r}}\times\left[  0,\beta\right[  )}^{2}\right)
\label{stima b2}%
\end{equation}
for some constant $c$ depending on $\lambda$ and $\overline{\Delta_{r}}%
\times\left[  0,\beta\right[  .$

\noindent Putting (\ref{stima b1}) and (\ref{stima b2}) in (\ref{B1_B2}) and
using (\ref{dis infinito}) we obtain (\ref{infinito locale}).

\ \ \ \ \ \ \ \ \ \ \ \ \ \ \ \ \ \ \ \ \ \ \ \ \ \ \ \ \ \ \ \ \ \ \ \ \ \ \ \ \ \ \ \ \ \ \ \ \ \ \ \ \ \ \ \ \ \ \ \ \ \ \ \ $\ \ \ \ \ \ \ \ \ \ \ \ \ \ \ \ \ \ \ \ \ \ \ \ \ \ \ \ \ \ \ \ \ \ \ \ \ \ \ \ \ \ \ \ \ \ \ \ \ \ \ \ \ \ \ \ \ \ \ \ \ \ \ \ \ \ \ \ \ \ \ \ \ \ \ \ \ \ \ \ \ \ \square
$\bigskip

\begin{remark}
\emph{In (\ref{traccia infinito}) the power }$2$ \emph{in the argument of the
exponential is sharp as the following example shows. In order to show that we
need to consider\ }$\Omega=\left\{  x\in\mathbb{R}^{N}:x_{N}<\omega\right\}
$\emph{ with }$\omega\in%
\mathbb{R}
$\emph{ and the function }$u_{\delta}(x)=u_{\delta}(x)=\left(  1-\log
\Phi(x_{N})\right)  ^{\delta}$\emph{ with} $0<\delta\leq\frac{1}{2}$ \emph{as
in Remark \ref{remark ottimalita infinito}} \emph{and argue as in Remark}
\emph{\ref{remark ottomalita traccia}}.
\end{remark}

\section{Trace operator}

\ \ \ \ In this section the "boundary values" or trace of functions in Sobolev
spaces are studied.

If $\Omega$ is a domain satisfying condition\emph{ \ref{domain}}, given a
smooth function $u\in C^{\infty}(\overline{\Omega})\subset W^{1,p}%
(\Omega,\gamma)$ we can define the restriction to the boundary $\left.
u\right\vert _{\partial\Omega}.$ It turn out that this restriction operator
can be extended from smooth functions to $W^{1,p}(\Omega,\gamma)$ giving a
linear continuous operator from $W^{1,p}(\Omega,\gamma)$ to $L^{p}%
(\partial\Omega,\gamma),$ the space of the measurable functions defined almost
everywhere on $\partial\Omega$ such that
\[
\int_{\partial\Omega}\left\vert u\right\vert ^{p}\varphi\text{ }%
d\mathcal{H}^{N-1}<+\infty.
\]
We stress that $L^{p}(\partial\Omega,\gamma)$ is a Banach space with respect
to the norm $\left\Vert u\right\Vert _{L^{p}(\partial\Omega,\gamma)}=\left(
\int_{\partial\Omega}\left\vert u\right\vert ^{p}\varphi\text{ }%
d\mathcal{H}^{N-1}\right)  ^{\frac{1}{p}}.$

Using the logarithmic Sobolev inequalities (\ref{dis lplogl}), there exists a
constant $C>0$ such that for every $u\in C^{\infty}(\overline{\Omega})$%
\begin{equation}
\left\Vert u\right\Vert _{L^{p}(\partial\Omega,\gamma)}\leq C\left\Vert
u\right\Vert _{W^{1,p}(\Omega,\gamma)}^{p}. \label{dis lp}%
\end{equation}
It follows that the operator
\begin{align*}
T  &  :C^{\infty}(\overline{\Omega})\rightarrow L^{p}(\partial\Omega,\gamma)\\
u  &  \rightarrow Tu=u/\partial\Omega
\end{align*}
is linear and continuous from $\left(  C^{\infty}(\overline{\Omega
}),\left\Vert {}\right\Vert _{W^{1,p}(\Omega,\gamma)}\right)  $ into $\left(
L^{p}(\partial\Omega,\gamma),\left\Vert {}\right\Vert _{L^{p}(\partial
\Omega,\gamma)}\right)  .$

By Hahn-Banach theorem and the density of $C^{\infty}(\overline{\Omega})$ in
$W^{1,p}(\Omega,\gamma)$ the operator can be extended to $W^{1,p}%
(\Omega,\gamma)$. This linear continuous operator from $W^{1,p}(\Omega
,\gamma)$ to $L^{p}(\partial\Omega,\gamma)$ is called trace operator of $u$ on
$\partial\Omega$. Then there exists a constant $C>0$ such that
\begin{equation}
\left\Vert Tu\right\Vert _{L^{p}(\partial\Omega,\gamma)}\leq C\left\Vert
u\right\Vert _{W^{1,p}(\Omega,\gamma)}\text{ \ \ for every }u\in
W^{1,p}(\Omega,\gamma), \label{dis traccia in W1p}%
\end{equation}
that implies that $W^{1,p}(\Omega,\gamma)$ is continuous imbedded in
$L^{p}(\partial\Omega,\gamma).$

Moreover the trace operator is compact for $1\leq p<+\infty.$ Indeed let
$\left\{  u_{n}\right\}  _{n\in%
\mathbb{N}
}$ be a bounded sequence in $W^{1,p}(\Omega,\gamma),$ we will prove the
existence of a Cauchy subsequence in $L^{p}(\partial\Omega,\gamma).$ By
Proposition \ref{proposizione compatezza}, there exists a Cauchy subsequence,
still denoted by $\left\{  u_{n}\right\}  _{n\in%
\mathbb{N}
},$ in $L^{p}\left(  \log L\right)  ^{\frac{1}{2p^{\prime}}}(\Omega,\gamma)$.
Moreover arguing as in the proof of the inequality (\ref{dis lplogl}) we have
\begin{align*}
\left\Vert Tu_{n}-Tu_{m}\right\Vert _{L^{p}(\partial\Omega,\gamma)}^{p}  &
\leq%
{\displaystyle\int_{\partial\Omega}}
\left\vert Tu_{n}-Tu_{m}\right\vert ^{p}\log^{\frac{p-1}{2}}(2+\left\vert
Tu_{n}-Tu_{m}\right\vert )\varphi\text{ }d\mathcal{H}^{N-1}\\
&  \leq c\left\Vert u_{n}-u_{m}\right\Vert _{L^{p}\left(  \log L\right)
^{\frac{1}{2}}(\Omega,\gamma)}^{p-1}\left\Vert \nabla\left(  u_{n}%
-u_{m}\right)  \right\Vert _{L^{p}(\Omega,\gamma)}\\
&  +c\left\Vert u_{n}-u_{m}\right\Vert _{L^{p}\left(  \log L\right)
^{\frac{p-1}{2}}(\Omega,\gamma)}^{p}\left\Vert x_{N}\right\Vert _{L^{\infty
}\left(  \log L\right)  ^{\frac{1}{2}}(\Omega,\gamma)},
\end{align*}
then $\left\{  u_{n}\right\}  _{n\in%
\mathbb{N}
}$ is a Cauchy sequence in $L^{p}(\partial\Omega,\gamma)$ too.

The norm of the trace operator is given by%

\begin{equation}
\underset{u\in W^{1,p}(\Omega,\gamma)-W_{0}^{1,p}(\Omega,\gamma)}{\inf}%
\frac{\left\Vert u\right\Vert _{W^{1,p}(\Omega,\gamma)}^{p}}{\left\Vert
Tu\right\Vert _{L^{p}(\partial\Omega,\gamma)}^{p}} \label{inf}%
\end{equation}
and this value is the best constant in the trace inequality
(\ref{dis traccia in W1p}). The trace operator is compact, therefore an easy
compactness arguments prove that there exist extremals in (\ref{inf}). These
extremals turn out to be the weak solution of
\begin{equation}
\left\{
\begin{array}
[c]{ll}%
-\operatorname{div}(\left\vert \nabla u\right\vert ^{p-2}\nabla u\varphi
)=\left\vert u\right\vert ^{p-2}u\varphi & \text{\emph{\mbox{in}} }\Omega\\
& \\
\left\vert \nabla u\right\vert ^{p-2}\frac{\partial u}{\partial\nu}%
=\lambda\left\vert u\right\vert ^{p-2}u\quad & \text{\emph{\mbox{on}}
}\partial\Omega,
\end{array}
\right.  \label{steklov 2}%
\end{equation}
where $\lambda$ is the first nontrivial eigenvalue.

\noindent When $p=2$ and $\Omega$ is a connected domain satisfying condition
\ref{domain}, using classical tools, compactness of the trace operator from
$W^{1,2}(\Omega,\gamma)$ to $L^{2}(\partial\Omega,\gamma)$\ and (\ref{poicare}%
) it follows that there exists an increasing sequence of eigenvalues of the
problem (\ref{steklov 2}) which tends to infinity and a Hilbertian basis of
eigenfunctions in $L^{2}(\Omega,\gamma).$

Moreover the continuity of the trace operator from $W^{1,2}(\Omega,\gamma)$ to
$L^{2}(\partial\Omega,\gamma)$\ and (\ref{poicare}) allow us to investigate
about the existence of a weak solution of the following semicoercive
nonhomogeneous Neumann problem
\[
\left\{
\begin{array}
[c]{ll}%
-\left(  u_{x_{i}}\varphi\right)  _{x_{i}}=f\varphi & \text{\emph{\mbox{in}}
}\Omega\\
& \\
\frac{\partial u}{\partial\nu}=g\quad & \text{\emph{\mbox{on}} }\partial
\Omega,
\end{array}
\right.
\]
where $\Omega$ is a connected domain satisfying condition \ref{domain}, $f\in
L^{2}(\log L)^{\frac{1}{2}}(\Omega,\gamma)$ and $g\in L^{2}(\partial
\Omega,\gamma).$ Indeed using classical tools (see e.g. \cite{libro buttazzo}
Theorem 6.2.5) we obtain that there exists a weak solution in $W^{1,2}%
(\Omega,\gamma)$ if and only if $\int_{\Omega}fd\gamma+\int_{\partial\Omega
}g\varphi$ $dH^{N-1}=0.$ In particular there exists a unique weak solution in
$X=\left\{  u\in W^{1,2}(\Omega,\gamma):\int_{\Omega}vd\gamma=0\right\}  $ by
Lax-Milgram theorem$.$

\section{Poincar\'{e} trace inequality}

\ \ \ \ Arguing as in Proposition \ref{prop poincare} (see Remark
\ref{remark sottospazi} too), we prove the following Poincar\'{e} type inequality.

\begin{proposition}
Let $\Omega$ be a connected domain satisfying condition \ref{domain} and
$1\leq p<+\infty.$ Then there exists a positive constant $C$, depending only
on $p$ and $\Omega$, such that
\begin{equation}
\left\Vert u\right\Vert _{L^{p}(\Omega,\gamma)}\leq C\left\Vert \nabla
u\right\Vert _{L^{p}(\Omega,\gamma)} \label{tipo poincare}%
\end{equation}
for any $u\in X=\left\{  u\in W^{1,p}(\Omega,\gamma):%
{\displaystyle\int_{\partial\Omega}}
u\varphi\text{ }d\mathcal{H}^{N-1}=0\right\}  .$
\end{proposition}

Using (\ref{dis lp}) and (\ref{tipo poincare}) we obtain

\begin{corollary}
Let $\Omega$ be a connected domain satisfying condition \ref{domain} and
$1\leq p<+\infty.$ Then there exists a positive constant $C$, depending only
on $p$ and $\Omega$, such that%
\begin{equation}
\left\Vert Tu\right\Vert _{L^{p}(\partial\Omega,\gamma)}\leq C\left\Vert
\nabla u\right\Vert _{L^{p}(\Omega,\gamma)} \label{poincare traccia}%
\end{equation}
for any $u\in X.$
\end{corollary}

\begin{remark}
\emph{(}Application to PDE\emph{)} \emph{Let consider the eigenvalue problem }%
\begin{equation}
\left\{
\begin{array}
[c]{ll}%
-\left(  u_{x_{i}}\varphi\right)  _{x_{i}}=0 & \text{\emph{\mbox{in}} }%
\Omega\\
& \\
\frac{\partial u}{\partial\nu}=\lambda u\quad & \text{\emph{\mbox{on}}
}\partial\Omega,
\end{array}
\right.  \label{steklov}%
\end{equation}
\emph{where }$\Omega$\emph{ is a connected domain satisfying condition
\ref{domain}. Arguing in a classical way using inequality
(\ref{poincare traccia}) and the compactness of the trace operator}$,$\emph{
it is easy to prove that} \emph{there exists an increasing sequence of
eigenvalues of the problem (\ref{steklov}) which tends to infinity. Moreover
for }$\lambda_{1}=0$\emph{ the corresponding eigenvalue function }%
$u_{1}=const\neq0$\emph{ and the first nontrivial eigenvalue }$\lambda_{2}%
$\emph{\ has the following characterization}%
\[
\lambda_{2}=\min\left\{  \frac{\left\Vert \nabla u\right\Vert _{L^{2}%
(\Omega,\gamma)}}{\left\Vert Tu\right\Vert _{L^{2}(\partial\Omega,\gamma)}%
},u\in W^{1,2}(\Omega,\gamma):\int_{\partial\Omega}u\varphi\text{
}d\mathcal{H}^{N-1}=0\right\}  .
\]

\end{remark}

\end{document}